\pdfoutput=1
\documentclass{article}
\usepackage[margin=25mm]{geometry}
\usepackage{amsmath}
\usepackage{amsfonts}
\usepackage{amssymb}
\usepackage{graphicx}
\usepackage[utf8]{inputenc}
\usepackage[T2A]{fontenc}
\usepackage[ukrainian,russian,english]{babel}
\usepackage[unicode, pdftex]{hyperref}
\usepackage{xcolor}
\usepackage[all]{xy}

\newtheorem{theorem}{Theorem}

\title{Typical one-parameter bifurcations of gradient flows with at most six singular points on the 2-sphere with holes}

\author{Svitlana Bilun, Maria Loseva, Olena Myshnova, Alexandr Prishlyak}

\begin{document}
\newcounter{contnumfig}
\setcounter{contnumfig}{0}

\maketitle
\begin{abstract}
   We describe all possible topological structures of typical one-parameter bifurcations of gradient flows  on the 2-sphere with holes in the case that the number of singular point of flows is at most six. To describe structures, we separatrix diagrams of flows. The saddle-node singularity is specified by selecting a separatrix in the diagram of the flow befor the bifurcation and the saddle connection is specified by a separatrix, which conect two saddles.
\end{abstract}
\section*{Introduction}

In this paper, we consider gradient flows on a sphere with holes. Since the function increases along each trajectory, the flow has no cycles and polycycles. In general position, a typical gradient flow is a Morse flow (Morse-Smale flow without closed trajectories). In typical one-parameter families of gradient flows, two types of bifurcations are possible: saddle-node and saddle connection. The corresponding vector fields at the time of the bifurcation are fields of codimension 1. In our case, they completely determine the topological type of the bifurcation.
To classify Morse flows, a separatrix diagram is often used, in which separetrices are trajectories that belonging to one-dimensiona stable or unstable manifolds.  We apply this approach to the classification of typical bifurcation.

Without loss of generality, we assume that under bifurcation (as the parameter increases), the number of singular points does not increase. The saddle-node bifurcation is defined by a separatrix, which is contructed to a point. We mark this separatrix on the diagram. A saddle connection bifurcation in the diagram corresponds to a separatrix, which conect  two saddles. 

We colar  stable separatrices in red, unstable separatrices in green and saddle connections in black.

 Kronrod \cite{Kronrod1950} and Reeb \cite{Reeb1946} construct topological invariants of functions oriented 2-maniofolds. It was generlized in \cite{lychak2009morse} for  non-orientable two-dimensional manifolds and in   \cite{Bolsinov2004, hladysh2017topology, hladysh2019simple} for manifolds with boundary, in \cite{prishlyak2002morse} for non-compact manifolds. 

In general, Morse flows are gradient flows of Morse functions. If we fix the value of functions in singular points the flow determinate the topological structure of the function \cite{lychak2009morse, Smale1961}. Therefore, Morse--Smale flows classification is closely related to the classification of the functions.

Topological classification of smooth function on closed 2-manifolds was obtained in \cite{bilun2023morseRP2, hladysh2019simple, hladysh2017topology,  prishlyak2002morse, prishlyak2000conjugacy,  prishlyak2007classification, lychak2009morse, prishlyak2002ms, prish2015top, prish1998sopr,  bilun2002closed,  Sharko1993}, on 2-manifolds with the boundary in \cite{hladysh2016functions, hladysh2019simple, hladysh2020deformations} and on closed 3- and 4-manifolds in  \cite{prishlyak1999equivalence, prishlyak2001conjugacy}.

In \cite{bilun2023gradient, Kybalko2018, Oshemkov1998, Peixoto1973, prishlyak1997graphs, prishlyak2020three, akchurin2022three, prishlyak2022topological, prishlyak2017morse,  kkp2013,  prish2002vek,  prishlyak2021flows,  prishlyak2020topology,   prishlyak2019optimal, prishlyak2022Boy}, 
the classifications of flows on closed 2- manifolds and 
\cite{bilun2023discrete, loseva2016topology, prishlyak2017morse, prishlyak2022topological, prishlyak2003sum, prishlyak2003topological, prishlyak1997graphs, prishlyak2019optimal} on manifolds with the boundary were obtained.
Topological properties of Morse-Smale flows on 3-manifolds was investigated in \cite{prish1998vek,  prish2001top, Prishlyak2002beh2, prishlyak2002ms,   prish2002vek, prishlyak2005complete, prishlyak2007complete, hatamian2020heegaard, bilun2022morse, bilun2022visualization}.


The purpose of this paper is to describe all possible topological structures of the Morse flows and typical bifurcations  with no more than six singular points (a saddle-node point we consider as two points) on a sphere with holes.



 
\section{Typical one-parameter bifurcations of gradient flows on a sphere with holes}

Typical vector fields on compact 2-manifolds are Morse-Smale fields. Morse fields or Morse-Smale gradient-like fields that do not contain closed trajectories are tipical among the gradient fields. They satisfy three properties:

1) singular points are nondegenerate;

2) there are no separatric connections;

3) $\alpha$-limiting ($\omega$-limiting) set of each trajectory is a singular point.

In typical one-parameter field families, one of these conditions is violated. Violation of the first condition leads to a saddle-node bifurcation, and the second to the appearance of a saddle connection. The third condition cannot be violated, because according to the Poincaré-Bendixon theorem, the $\alpha$-limit ($\omega$-limit) set of every trajectory on the sphere is either a singular point, or a cycle, or a polycycle. Since gradient fields have no cycles and polycycles, this set is a singular point.

\subsection{Internal bifurcations}

According to the theory of bifurcations, there are only two typical bifurcations of gradient flows: a saddle-node and a saddle connection.

\subsubsection{Saddle-node bifurcation}
\begin{figure}[ht!]
\center{\includegraphics[width=0.95\linewidth]{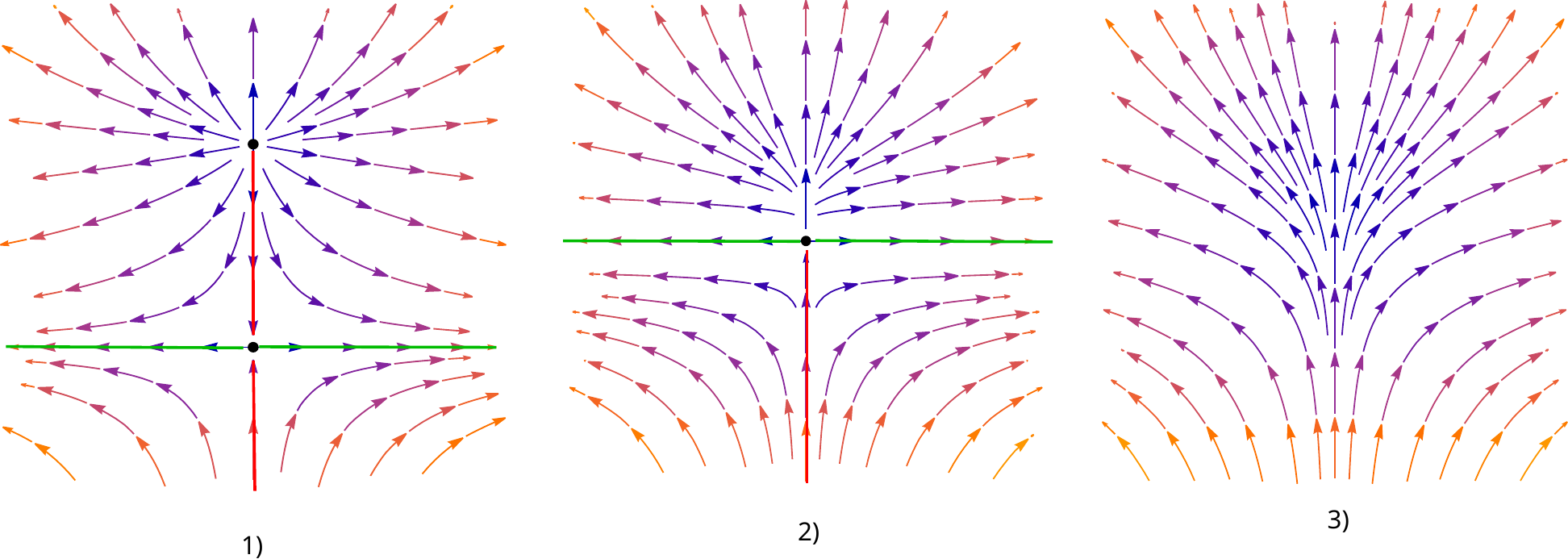}
}
\caption{SN -- saddle-node bifurcation
}
\label{bifsn}
\end{figure}

The saddle-node bifurcation, when the node is the source, is shown in Fig. \ref{bifsn}.

It can be described by the equation $V(x,y,a)=\{x,y^2+a\}$ if the node is the source and the equation $V(x,y,a)=\{-x,- y^2a\}$ if the node is a sink. Here, $a$ is a parameter. If $a<0$ we get the flow before the bifurcation, if $a>0$ we get the flow after the bifurcation, and if $a=0$ we get the flow at the moment of the bifurcation (the flow of codimensionality 1). The saddle-sink bifurcation can be obtained from the saddle-source by changing the direction of all trajectories, that is, by taking the reversed flow. Note that the saddle-node bifurcation corresponds to a typical one-parameter
catastrophe (deformation) of functions: $f(x,y,a)=x^2+y^3+ay$.

Consider a stable manifold (red separatrix in Fig. 1 
for a singular saddle-source point. 
In order to determine the saddle-node bifurcation it is necessary to select a separatrix on the seperatrix diagram.

Note that the bifurcation of the saddle flow in Morse's theory corresponds to the operation of reducing points or reducing complementary handles.

\subsubsection{Saddle connection}
\begin{figure}[ht!]
\center{\includegraphics[width=0.95\linewidth]{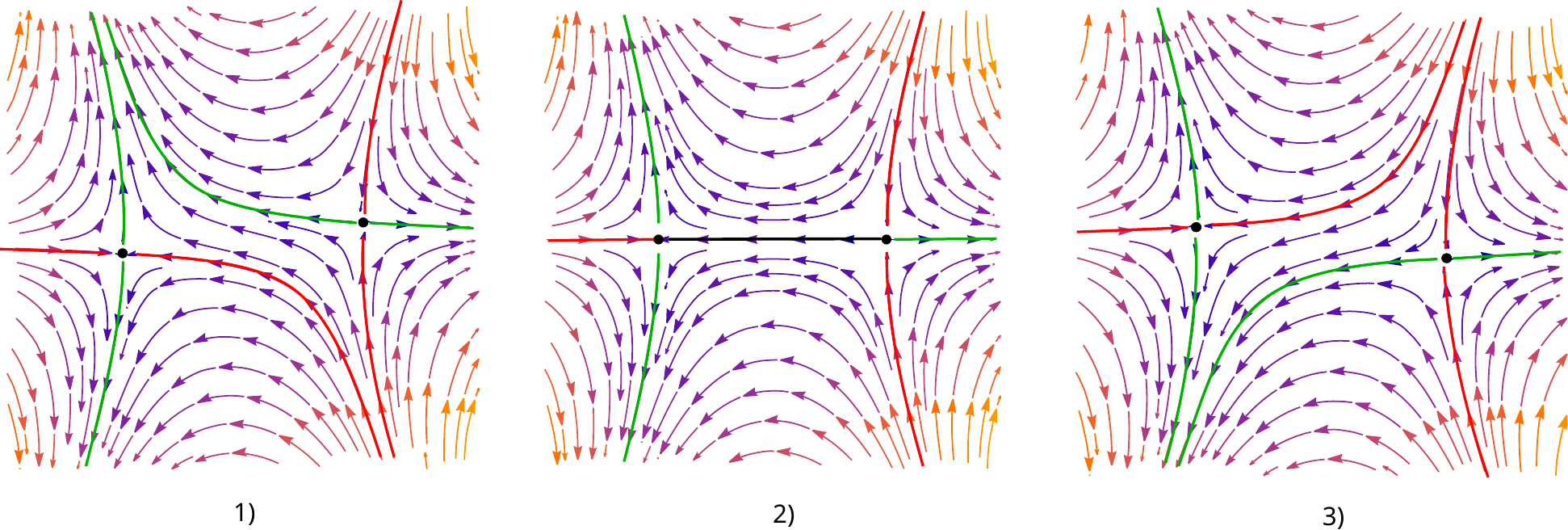}
}
\caption{SC -- saddle connection
}
\label{bifsc}
\end{figure}

The bifurcation of the saddle connection is shown in 
Fig. \ref{bifsc}. 
It can be described by the equation $V(x,y,a)=\{x^2-y^2-1,-2xy+a\}$.

We describe one of the possible situations when it appears in typical one-parameter families of gradient fields. Let $p,q$ be saddle critical points of the function $f$, $f(p)<f(q)$, $[\frac{f(p)+f(q)}{2}-\varepsilon, \frac{f(p)+f(q)}{2}+\varepsilon]$ does not contain critical values of $f$, $L$ of the line component of level $f^(-1)(\frac{f(p) +f(q)}{2})$, $u$ is the separatrix of $p$, $v$ is the separatrix of $q$, which intersect $L$.
Then according to Morse theory, $f^{-1}([\frac{f(p)+f(q)}{2}-\varepsilon, \frac{f(p)+f(q)}{2} +\varepsilon])$ is homeomorphic to the cylinder $S^1 \times [0,1]$ and in the coordinates $(s,t)$ ($s$--polar angle on the circle) on it, the trajectories are given by segments $s= \text{const}$. Consider the twisting of a cylinder given by the formula $(s,t) \to (s+a,t)$ ($a$-- deformation parameter). At the same time, the trajectories turn into helical lines. With a continuous change of the parameter $a$, the point of intersection of the trajectory $u$ with the circle $frac{f(p)+f(q)}{2}+\varepsilon$ continuously rotate along it. The intersection point of the trajectory $v$ with this circle remain unchanged. By the intermediate value theorem, there is a moment in time $a$ at which these points coincide. It is at this point in time that the bifurcation of the saddle connection occur.
Note that the construction described above can be implemented by changing the Riemannian metric (changing angles considered perpendicular), without changing the function itself.
The saddle connection will be indicated by black separatrix in the separatrix diagram.  

In Morse's theory, bifurcation of the saddle connection corresponds to the operation of adding handles (sliding one handle over another handle of the same index).

\subsection{Bifurcations of singular points on the boundary}

With such bifurcations, two cases are possible: 1) the topological structure of the flow restriction to the boundary changes, 2) this structure does not change.

Using the Poincaré-Hopf theorem for the doubling of the flow, we can conclude that the case when one point on the boundary disappears or two points on the boundary are replaced by one is not possible, because if this is the case, then the parity of the sum of the indices of the flow on the doubling is violated. Therefore, the only way to change the structure of the flow at the boundary in a typical family is the coalescence of two points at the boundary, which after the bifurcation disappear or are replaced by one internal singular point.

In the second case, as for the internal saddle-node bifurcation, the separatrix connecting them is compressed to a point. At the same time, one of the two points before the bifurcation (saddle or node) lies on the boundary. Both points cannot lie on the boundary, because then the structure of the manifold at the point obtained after the compression of the separatrix will be disturbed after compression.

Bifurcations that do not change the flow at the boundary also include bifurcations of a saddle connection, when one or both saddles lie on the boundary. The only thing that should be taken into account here is that these saddles cannot be adjacent points on the boundary (there are other singular points between them). Otherwise, the trajectory at the boundary between them, together with the separatrix connection at the moment of bifurcation, form a polycycle of length 2, which is not possible for gradient flows.

Depending on the types of points that stick together, different options for bifurcations are possible.

\begin{figure}[ht!]
{\center{\includegraphics[width=0.9\linewidth]{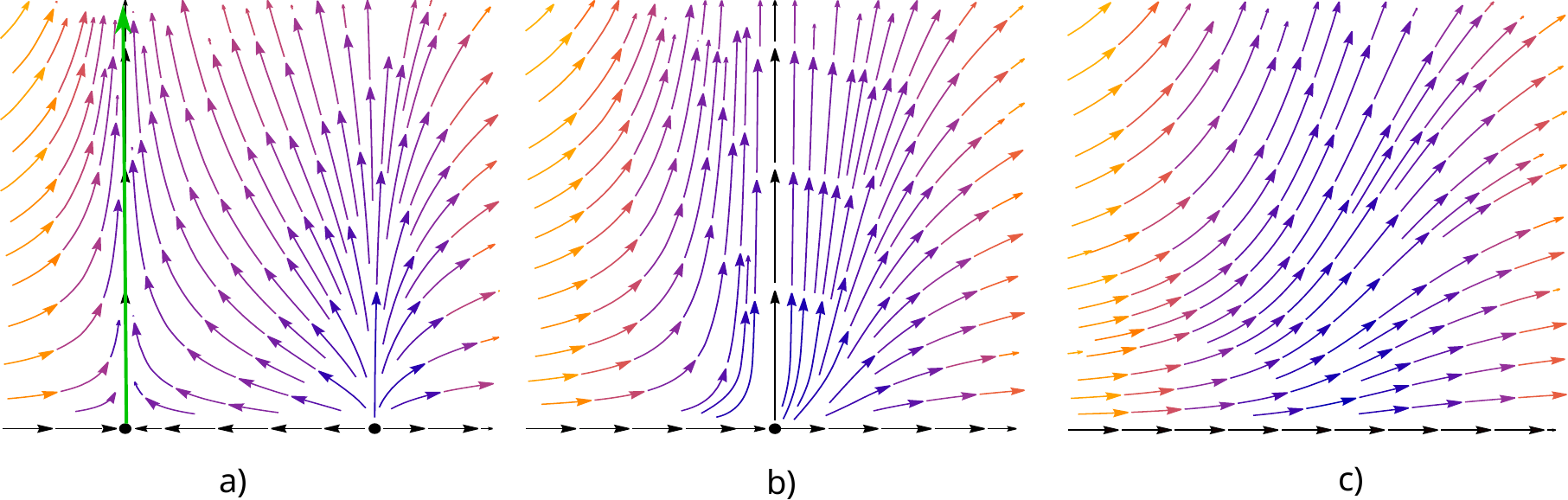}}}
\caption{ BSN -- saddle-node (source) bifurcation at the surface boundary }
\label{s-d}
\end{figure}

1) In the first case, the sourse and saddle are glueded together at a point.
In Fig. \ref{s-d} a) we show the flow before the bifurcation ($a=-1$), in Fig. \ref{s-d} b) -- flow at the moment of bifurcation ($a=0$), in Fig. \ref{s-d} c) -- flow after bifurcation ($a=1$).

\begin{figure}[ht!]
{\center{\includegraphics[width=0.9\linewidth]{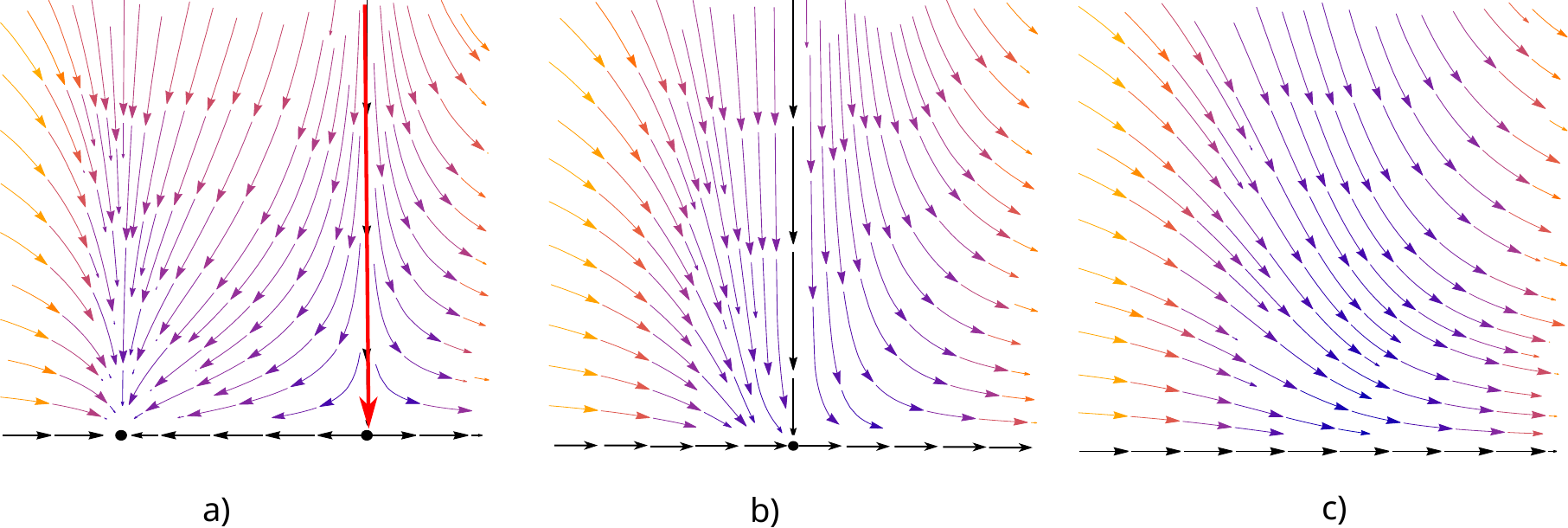}}}
\caption{ BSN -- saddle-node bifurcation (sink) at the surface boundary }
\label{s-k}
\end{figure}

2) In the second case, the saddle and the sink merge into a point, which then disappears (Fig.\ref{s-k}).

If one of the two saddle-node bifurcation points is internal, and the other lies on the boundary, then we have two types of semi-boundary saddle-node bifurcations: at the moment of bifurcation, a saddle (HS) or a node (HN).

\begin{figure}[ht!]
{\center{\includegraphics[width=0.9\linewidth]{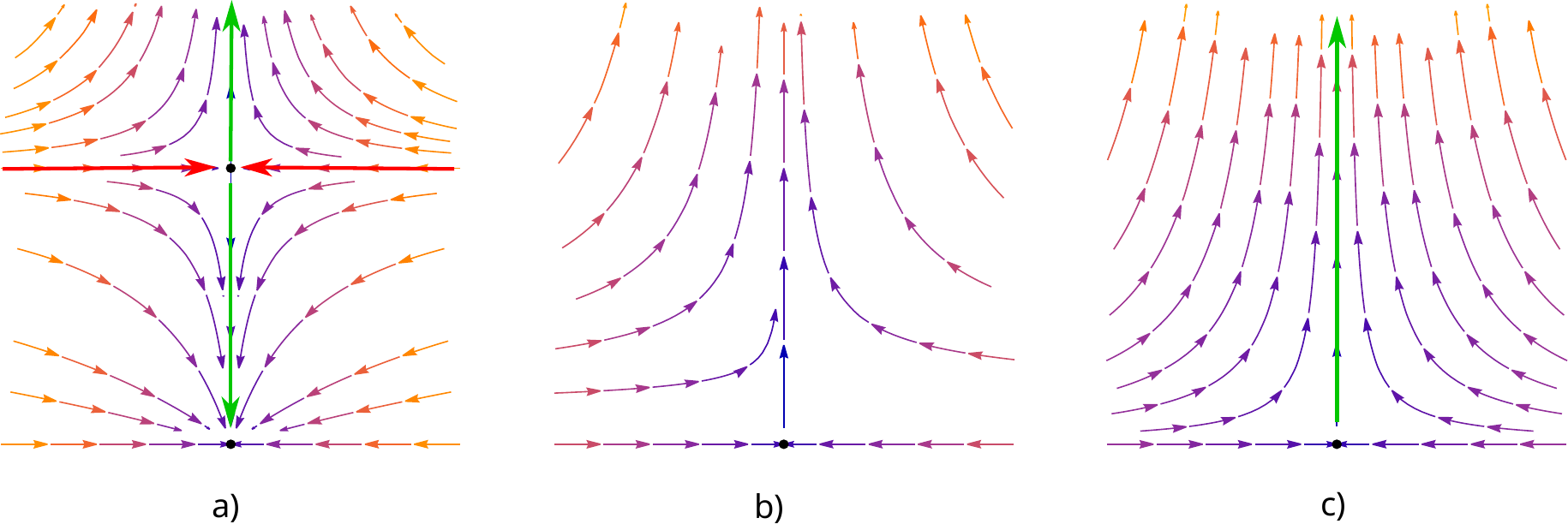}}}
\caption{ HS -- semi-boundary saddle bifurcation   }
\label{HS}
\end{figure}

\newpage

\begin{figure}[ht!]
\center{\includegraphics[width=0.9\linewidth]{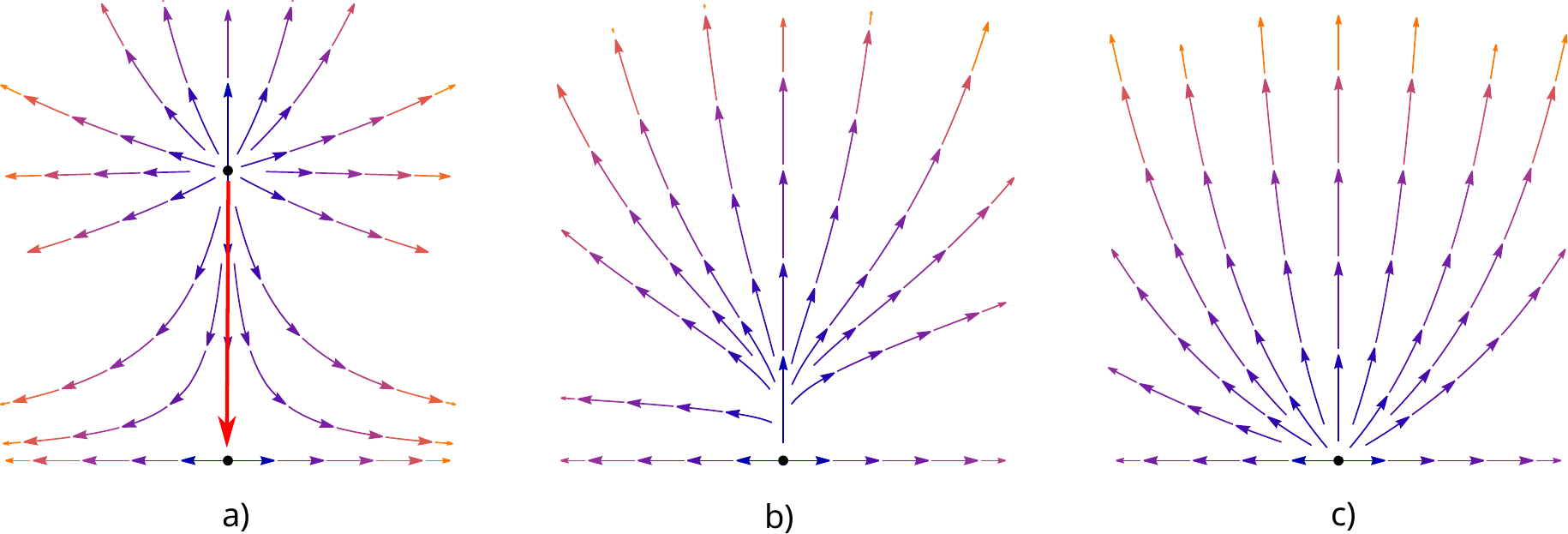}}
\label{HN}
\caption{HN -- semi-boundary node bifurcation }
\end{figure}

\begin{figure}[ht!]
\center{\includegraphics[width=0.9\linewidth]{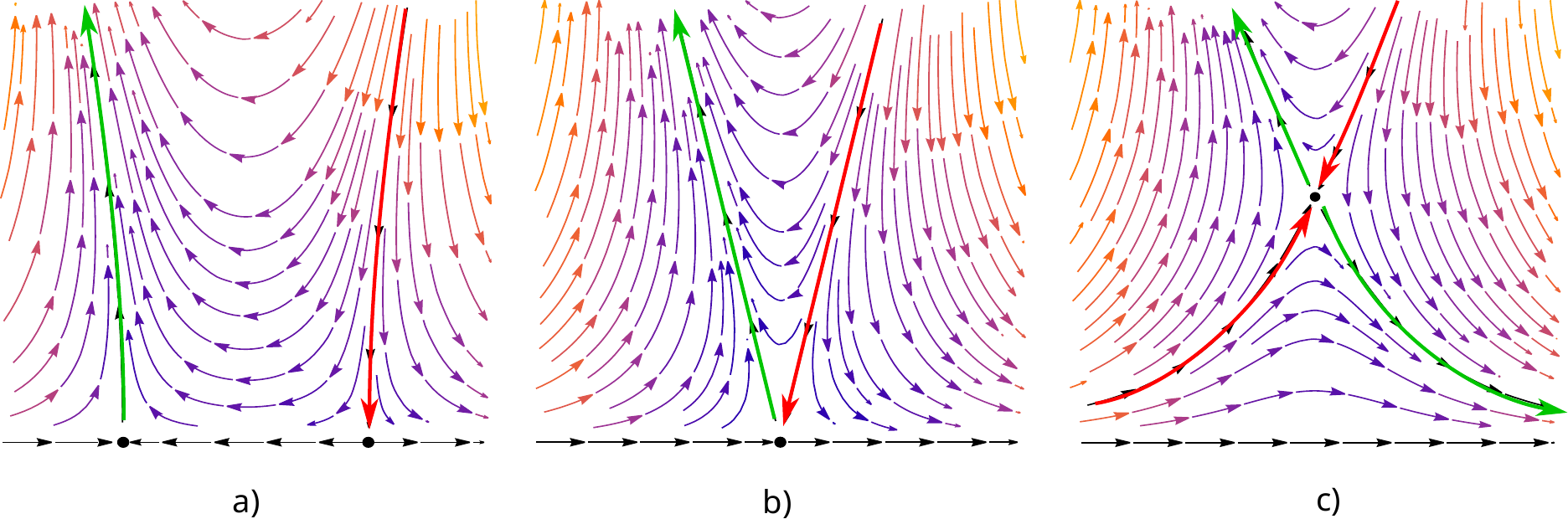}}
\caption{BDS -- double saddle bifurcation at the boundary }
\label{s-s}
\end{figure}

In addition, the following options are possible:

3) both points that stick together are saddles:
  $\{ x^2-y^2+a, -2xy \}$ (Fig. \ref{s-s}),

\begin{figure}[ht!]
\center{\includegraphics[width=0.9\linewidth]{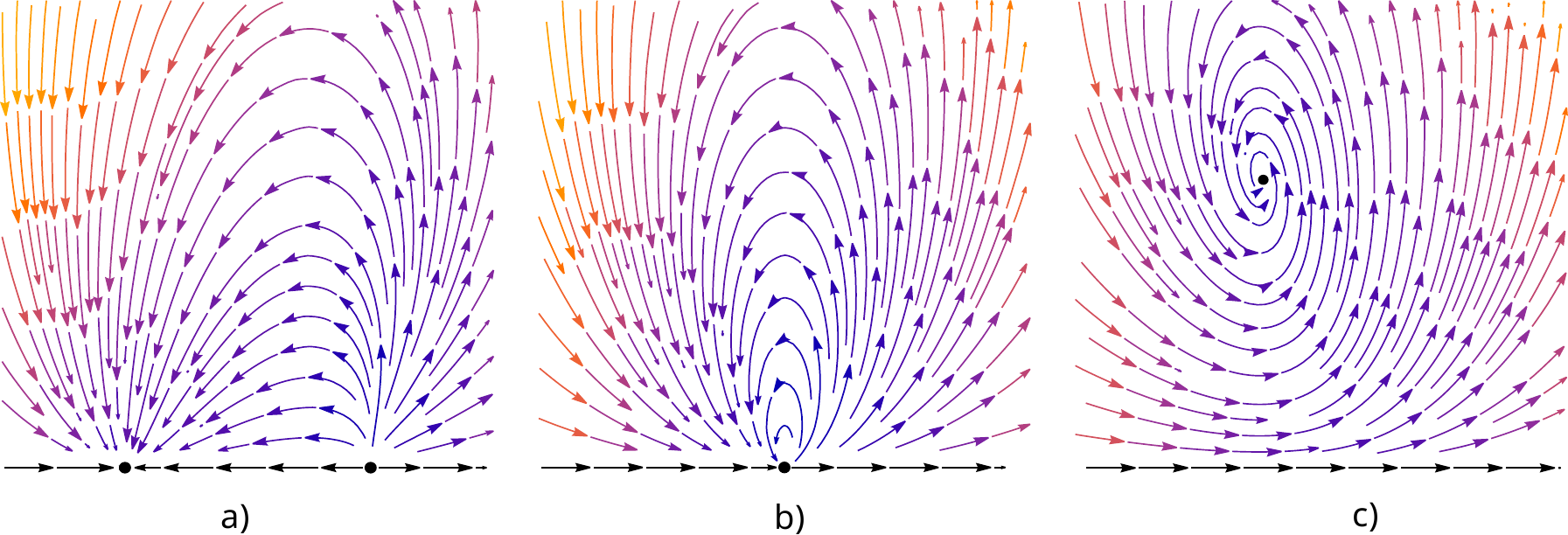}}
\caption{ source-sink bifurcation -- sink}
\label{d-k}
\end{figure}

4) source and sink stick together, an internal sink appears
$\{ x^2-y^2+a, 2xy+ay \}$ \ (Fig. \ref{d-k})

\newpage

\begin{figure}[ht!]
\center{\includegraphics[width=0.9\linewidth]{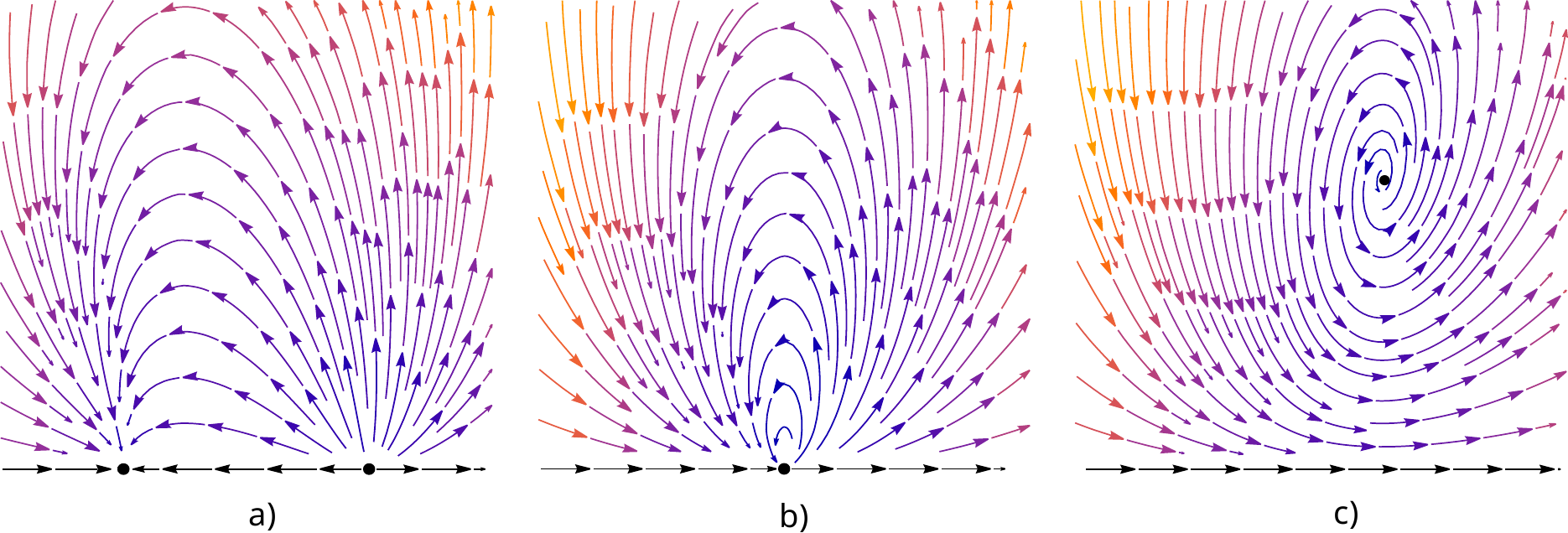}}
\caption{ source-sink bifurcation -- source}
\label{d-k2}
\end{figure}

5) source and sink stick together, an internal source appears
$\{ x^2-y^2+a, 2xy - ay \}$ \ (Fig. \ref{d-k2})

In the last two variants, an internal singular point appears, since when doubling the vector field on the sphere according to the Poincaré-Hopf theorem, the sum of the Poincaré indices does not change during the bifurcation, then two saddles on the boundary (index -1) will give rise to one saddle in each of the half-planes, and source and sink (index +1) source or sink in each half-plane. Therefore, only the saddle-node bifurcation does not lead to the appearance of internal singular points. Note that all these bifurcations at the moment $a=0$ define a flow with fixed points on the boundary.

As for manifolds without a boundary, in addition to local bifurcations, in typical one-parameter families of vector fields, only one type of global bifurcations is possible - a saddle connection

We will show that the described local bifurcations occur in typical one-parametric families of flows. Since the saddle-node bifurcation occurs in typical one-parameter families of flows on a plane, when it is restricted to a half-plane, it will give a typical one-parameter family on a half-plane.

Consider the situation of the merger of two saddles. The flow on the line $y=0$ has a sink at $x=-1$ and a source at $x=1$.

It continues to the following saddles: $\{ x^2-1,y\}$ in the vicinity of $(-1,0)$ and $\{ x^2-1,-y\}$ in the vicinity of $(1, 0)$.

The simplest function that makes a continuous transition from 1 (at $x=-1$) to -1
(at $x=1$) there is a function $ - x$.
Therefore, the flow formula will be written as $\{ x^2-1,-xy \}$.
If we now consider the bifurcation $\{ x^2-a, -xy \}$, then at $a=0$ a whole straight line $x=0$ will arise from zeros of the flow, which is not possible according to our assumptions. To prevent this from happening, let's add to the first coordinate a function that rotates to zero on the line $y=0$, the simplest such functions are $y$ and $-y$, and for symmetric fields with respect to the horizontal axis, $y^2$ and $ -y^2$. Direct inspection shows that the first functions lead to an interior singular point on the upper half-plane. So, the simplest bifurcation has the form $\{ x^2-y^2+a, -2xy \}$ (for convenience, we stretched the field twice as much along the y axis). In the general case, we can also change the second coordinate of the vector field, but considering that this coordinate is zero at $y=0$, the changes will be written as $\{ x^2-y^2+a, -2xy +cay \}$ for of some constant $c$. But for all $c$ we will get the same bifurcations.

When the source and the sink stick together, we will perform the same considerations and obtain a bifurcation of a typical one-parametric family according to the formula $\{x^2-y^2+a, 2xy +cay \}$. Positive and negative values of $c$ give two different bifurcations (described above), and with $c=0$ we get a non-typical family of fields.

\begin{theorem} For typical one-parametric families of flows without closed trajectories and with fixed points on the boundary, all bifurcation points are isolated.
\end{theorem}
\textbf{Proof.} In the flat case, the bifurcation point can be the point of accumulation of bifurcations, if the cycle breaks up at the bifurcation - a closed curve or a homoclinic trajectory. Since we assume that the fields do not have closed trajectories, the first case is not possible. In the second case, the homo clinical trajectory begins and ends in the saddle. This is also not possible for a saddle on the boundary, for example, if the saddle has the form $\{ x,-y \}, y \ge 0$, then the trajectories leaving the saddle belong to the boundary and cannot enter this saddle, since it includes only one trajectory, which is internal. Therefore, both cases of accumulation of bifurcation points are not possible. The theorem is proved.

Therefore, for each bifurcation, all flows at $a<0$, as well as all flows at $a>0$, are topologically equivalent to each other. Without limiting the generality, we will assume that the number of flow zeros at a<0 is not less than at $a>0$. The stream at $a=-1$ will be called the initial stream.

\begin{theorem} On a set of flows with fixed points on the boundary and without closed trajectories, the bifurcation in a typical family is given either by the initial flow and a compressible trajectory for local bifurcations on the boundary or by the flow at the moment of bifurcation in the case of a saddle connection. In case of merging of source and sink after bifurcation, a flow with internal zero occurs. To determine such a bifurcation, it is necessary to additionally specify whether a source or a sink is formed.
\end{theorem}

\begin{theorem} The first three bifurcations can be realized in the space of gradient flows, and the last two do not allow such realizations.
\end{theorem}
\textbf{Proof.} The one-parameter family of functions (catastrophe) corresponding to the first bifurcation has the form $f(x,y,a)=x^3/3+ax+y^2/2$, and for the second - - $f(x,y,a)=x^3/3+ax-y^2/2$. For the third bifurcation, the family of functions is given by the formula $f(x,y,a)=x^3/3+ax-xy^2$.

For the third and fourth bifurcations at $a=0$, the flow consists of loops, which is not possible for the gradient flow, since the value of the function increases along each trajectory of the gradient flow.

Moreover, on the set of gradient flows there is no bifurcation between the initial and final values ($a=-1$ and $a=1$) of the indicated bifurcations. Indeed, suppose that such a bifurcation exists. Let us continue the flow at $a=-1$ to the optimal flow of the Morse flow on a compact surface as in Losev-Pryshlyak. This flow has a single sink and a single source. The corresponding Morse function will have one local minimum and one local maximum. If we assume that there is a bifurcation of the general position, then as a result we will get a function without a local minimum or a local maximum on a compact surface, which is not possible.

Therefore, the following types of gradient bifurcations are possible on spheres with holes:

SN -- internal saddle node;

SC -- internal saddle connection;

BSN - border saddle node;

BDS -- limit double saddle;

HN -- semi-boundary saddle node (node);

HS -- semi-boundary saddle-node (saddle);

HSC -- semimarginal saddle connection;

BSC -- a bundle of saddles on the border.

In the case of saddle-node bifurcations, such a bifurcation is given by a separatrix diagram to the bifurcation, on which the trajectory (separatrix) between the saddle and the node is highlighted, which is compressed to a point. To specify the bifurcation of the saddle connection, a separatrix diagram at the moment of bifurcation is sufficient.

\section{The structure of typical flows and bifurcations on a 2-disk}

In Fig. 10, we show all possible (with accuracy to homeomorphism) separatrix flow diagrams
Morse on 2-disc with no more than 4 singular points.

\begin{figure}[ht!]
\center{\includegraphics[width=0.95\linewidth]{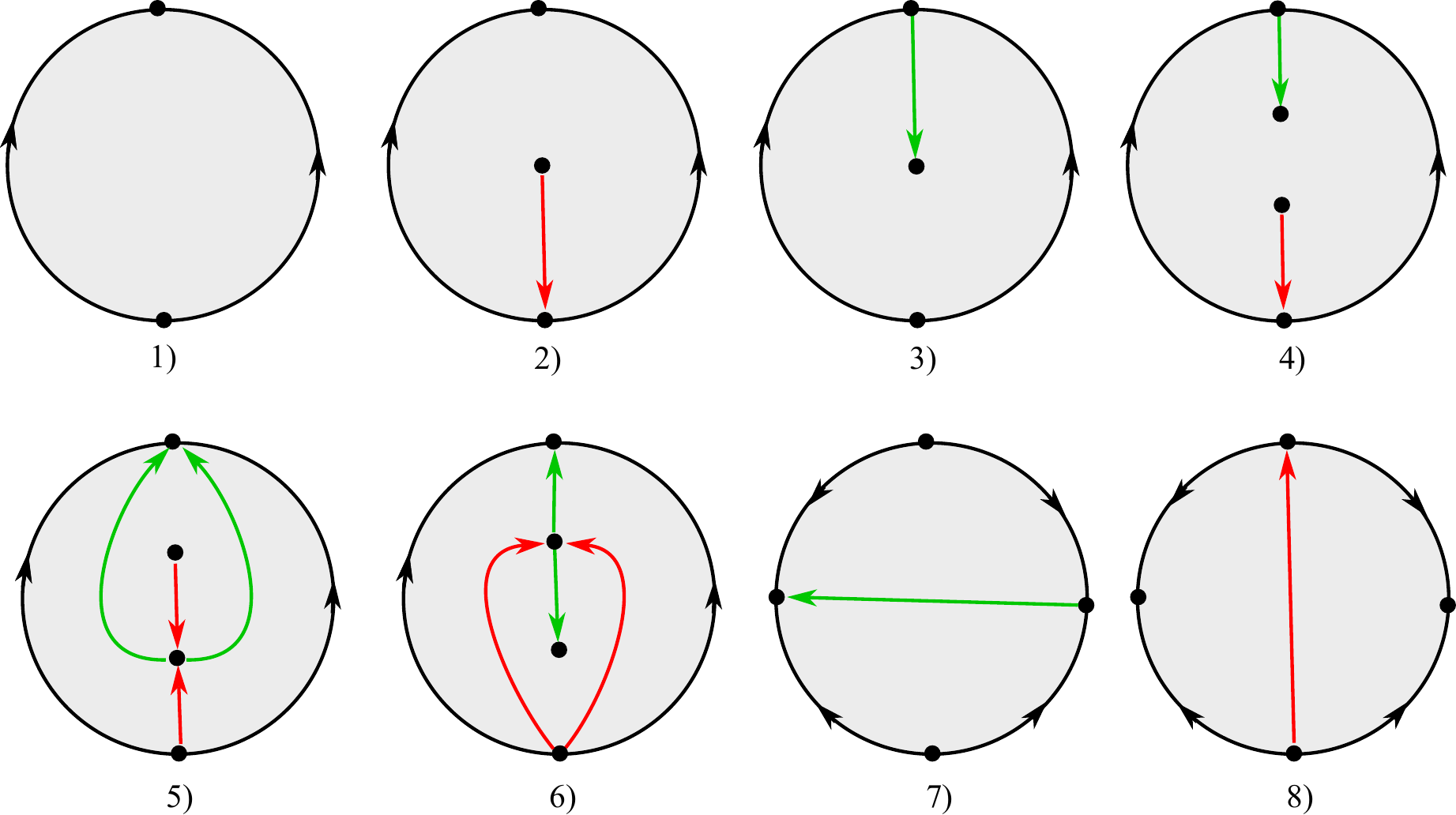}
}
{
\caption{Morse flows on 2-disc with no more than 4 singular points}
}
\label{D2-24}
\end{figure}
At the same time, diagrams 2) and 3), 5) and 6), as well as 7) and 8) are mutually inverted. In the following figures, we will depict only one of such diagrams (the one with the largest number of leaks). As can be seen from the figures, only diagram 4) contains two saddles. Therefore, bifurcation of saddle connections is impossible for flows with one saddle, and it is also not possible for flow 4) because the saddles are adjacent on the boundary.

For saddle-node bifurcations, it is only necessary to note how many different separatrixes and limit trajectories connecting a saddle and a node exist (with homeomorphism accuracy) on each diagram.

Bifurcations are not possible for the first diagram. For diagrams 2) and 3) we have one bifurcation (HN) each, for 4) two such bifurcations (HN), for 5) and 6) one SN each, two HS each. Boundary (both points on the boundary) for the indicated diagrams are not possible, because if one of the two arcs at the moment of bifurcation contracts to a point, then the other will form a loop, which is not possible for gradient flows.

For diagrams 7) and 8) we have one BSN bifurcation each (arcs connecting the saddle and the node on the boundary).

Next, we consider flows with five singular points (Fig. 11).

\begin{figure}[ht!] 
\center{\includegraphics[width=0.95\linewidth]{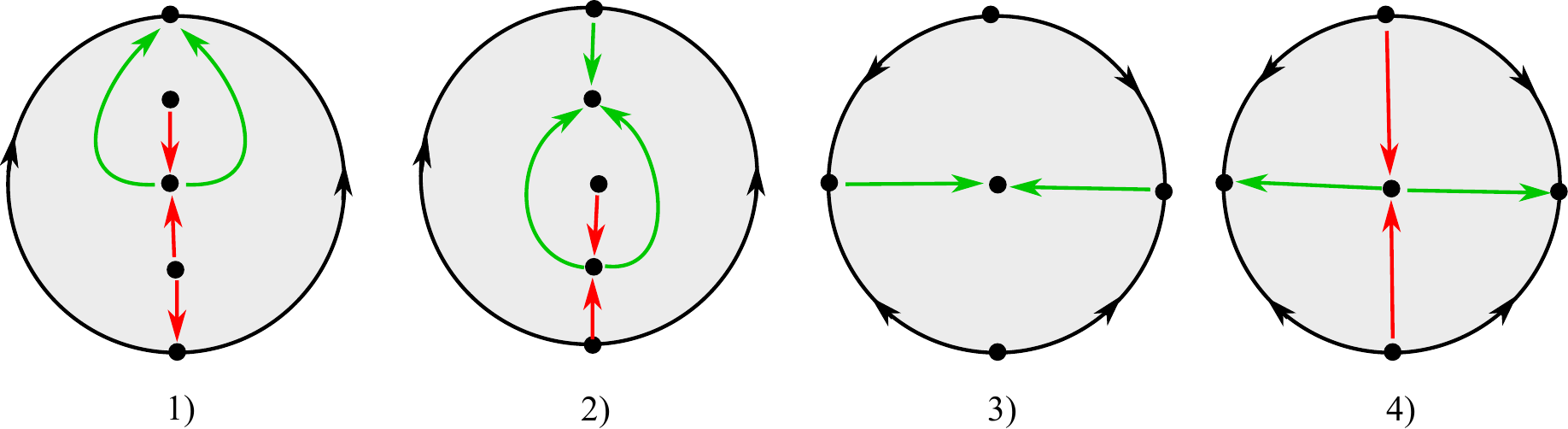}
}
{
\caption{Morse flows with 5 singu;ar point on the 2-disk}
}
\label{D2-5}
\end{figure}

The following saddle-node bifurcations are possible for these flows:

1) 2 SN, 1 HN, 1 HS;

2) 2 SN, 1 HN, 1 HS;

3) 1 BSN, 1 HN;

4) 2 HS.

Note that only 4) of the considered diagrams will turn into itself when the flow is reversed. Therefore, the total number of bifurcations in the first three cases is multiplied by 2, and for the last one we leave it unchanged.

In fig. 12 shows all possible separatrix flow diagrams with 6 singular points.
\begin{figure}[ht!]
\center{\includegraphics[width=0.95\linewidth]{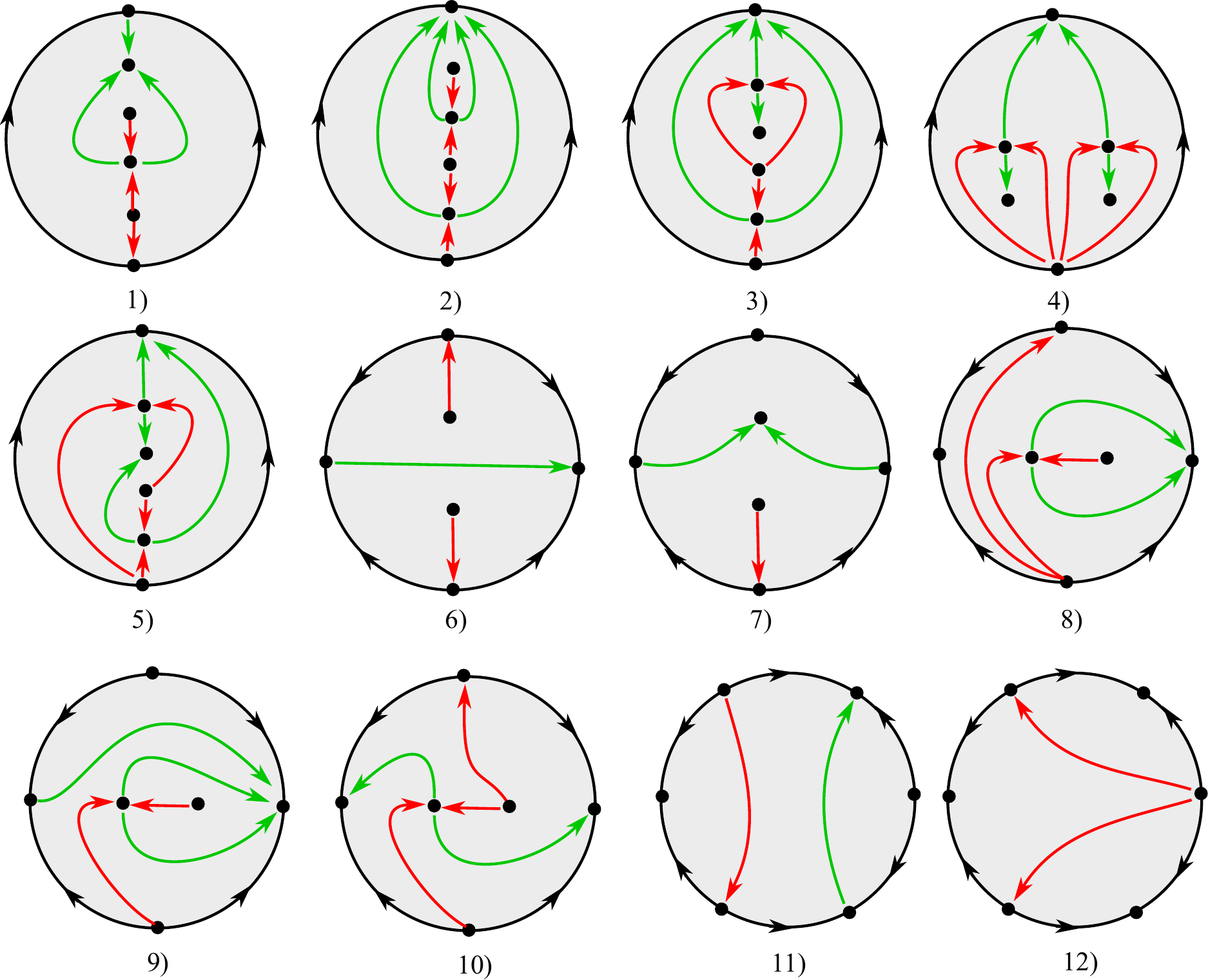}
}
{
\caption{Morse streams on 2-disk with 6 singular points }
}
\label{D2-6}
\end{figure}
Depending on the number of separatrixes on each of the diagrams, we get the following possible bifurcations for each of them:

1) 3 SN, 2 HN;

2) 3 SN, 3 HS;

3) 3 SN, 3 HS;

4) 1 SN, 3 HS;

5) 4 SN, 4 HS;

6) 1 HN, 1 BDS;

7) 2 HN, 1 BDS, 1 BSN;

8) 1 SN, 2 BSN, 2 HS;

9) 1 SN, 3 HS, 2 BSN;

10) 1 SN, 3 HS, 1 HN, 2 BSN;

11) 1 BDS, 2 BSN;

12) 2 BSN.

When changing the direction of the flow, diagrams 5) and 11) will not change, therefore, in the general calculation of the number of bifurcations, the above values will not change. For the rest of the charts, these values will be doubled.

All possible separatrix connections on a two-dimensional disc with no more than 6 singular points are shown in fig. 13.

\begin{figure}[ht!]
\center{\includegraphics[width=0.95\linewidth]{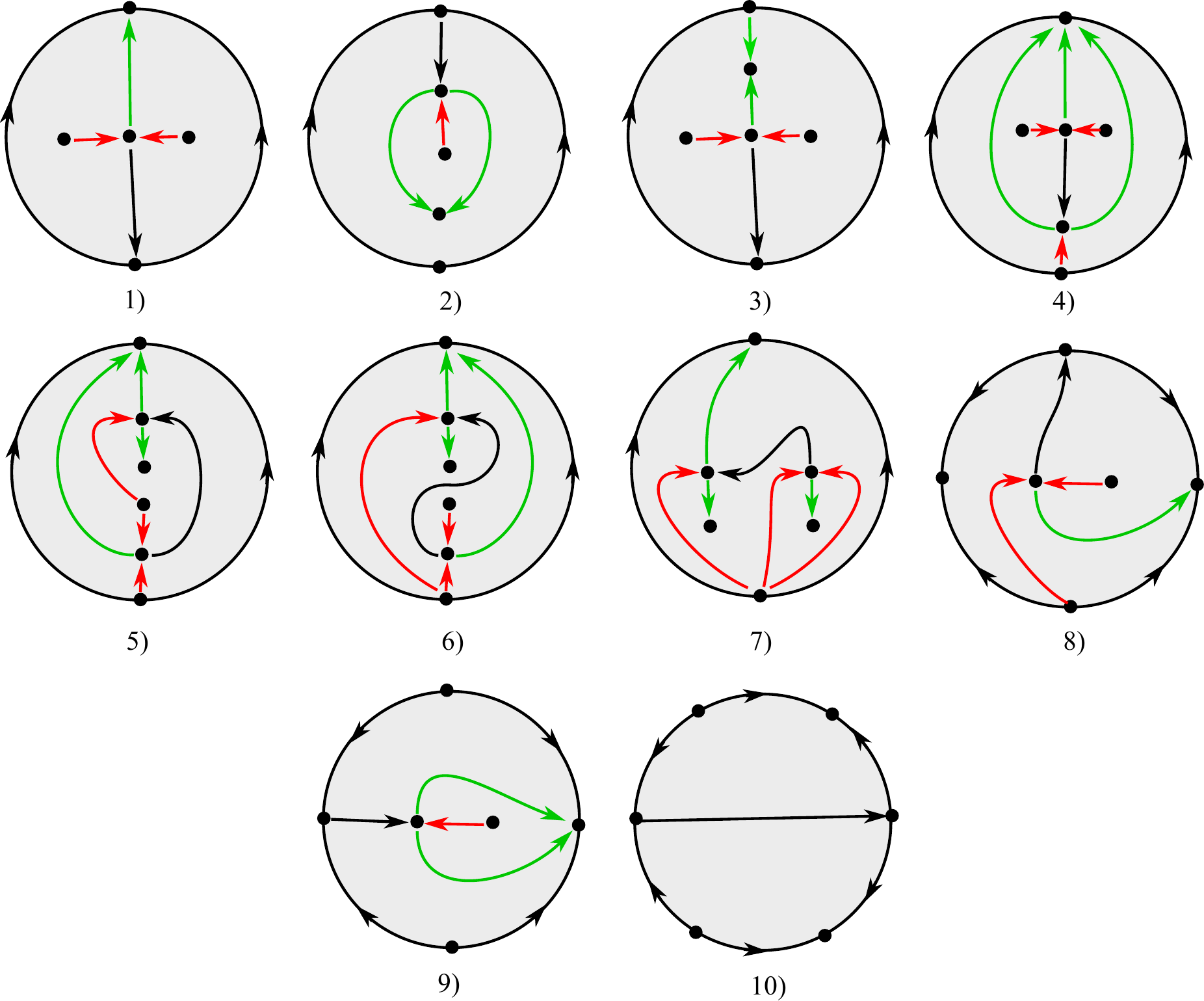}
}
{
\caption{Separatrix connections on a 2-disc with no more than 6 singular points}
}
\label{D2-sc}
\end{figure}

We went through all the possible options, and therefore it is fair

\begin{theorem}
The following possible structures of typical one-parameter gradient saddle-node bifurcations exist on a two-dimensional disk:
\begin{itemize}
\item
with two singular points at the moment of bifurcation - two bifurcations of HN;
\item
with three singular points at the moment of bifurcation - two bifurcations HN, two SN, four HS and two BSN;
\item
with four singular points at the moment of bifurcation -8 SN, 6 HN, 2 BSN, 8 HS;
\item
with five singular points at the moment of bifurcation - 30 SN, 12 HN, 38 HS, 22 BSN, 5 BDS.
\end{itemize}

There are the following bifurcations of the saddle connection:
\begin{itemize}
\item
with five singular points - four HSC bifurcations ( 1)$\times$2 and 2)$\times$2 in fig. 13);
\item
seven SC (4)$\times$2,5)$\times$2,6),7)$\times$2), six HSC (3)$\times$2, 8$\times$2, 9)$\times$2 ) and two BSC (10)$\times$2).
\end{itemize}

\end{theorem}

\section{Structure of typical flows and bifurcations on a cylinder}

Since there must be at least two singular points on each component of the boundary, then if the total number of such points is 4, then all of them lie on the boundary.

If the total number of singular points is 5, then one of them is internal, and two singular points lie on each component of the boundary (on each component of the boundary there is an even number of points, because when limiting the flow to the boundary, sources and sinks alternate).

\begin{figure}[ht!]
\center{\includegraphics[width=0.95\linewidth]{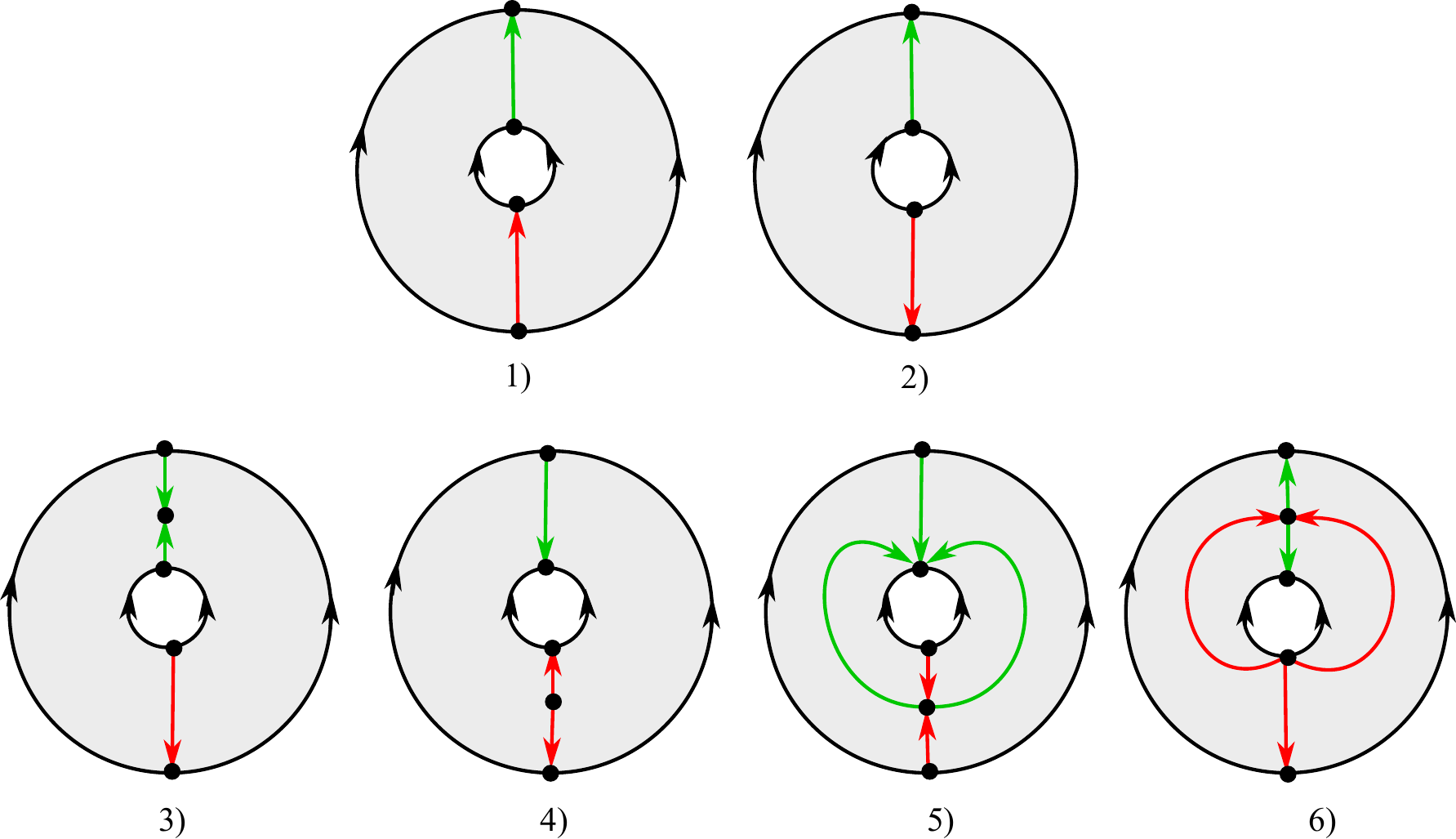}
}
{
\caption{Morse flows on a cylinder with 4 and 5 singular points}
}
\label{SI-45}
\end{figure}

All possible separatrix diagrams of Morse flows with 4 and 5 singular points are shown in Fig. 14. For the first two diagrams, saddle-node bifurcations are not possible, since imposing an arbitrary separatrix on them will change the topological type of the surface. The following saddle-node bifurcations are possible on other diagrams:

3) 2 HN;

4) 2 HN;

5) 3 HN;

6) 3 HN.

In fig. 15. shows all the separatrix diagrams of Morse flows on a cylinder with six singular points. For each pair of flow and its inverse, we have depicted only one diagram.

\begin{figure}[ht!]
\center{\includegraphics[width=0.95\linewidth]{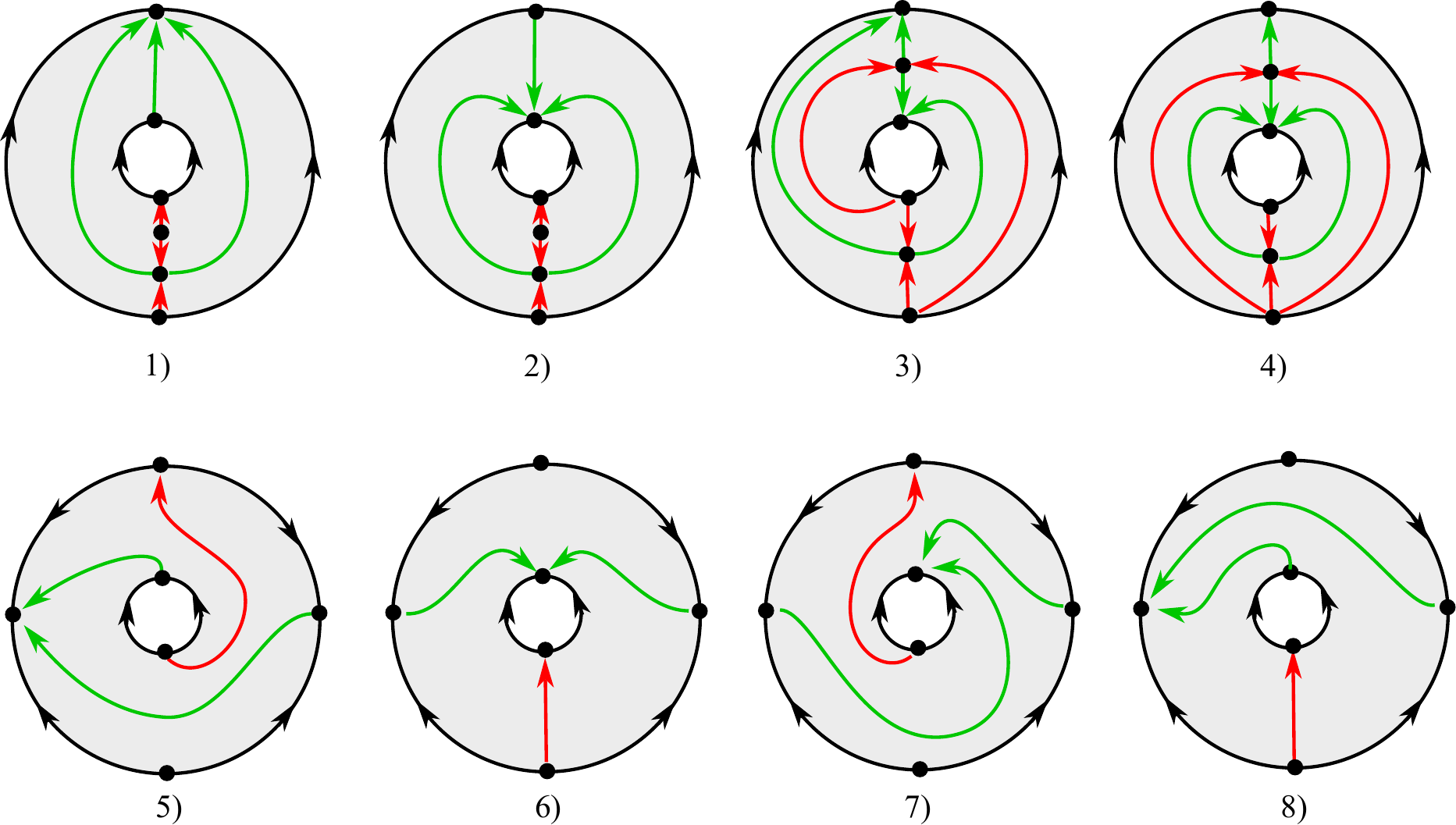}
}
{
\caption{Morse flows on a cylinder with 6 singular points}
}
\label{SI-6}
\end{figure}

The following bifurcations are possible for them:

1) 1 SN, 2 HS, 1HN;

2) 1 SN, 2 HS, 1HN;

3) 4 HS;

4) 6 HS;

5) 2 BSN, 1 BDS;

6) 2 BSN;

7) 1 BSN, 1 BDS;

8) 2 BSN, 1 BDS;

Note that only the third and fourth diagrams will turn into themselves when the flow orientation is changed.

All possible saddle connections of threads with 4 and 5 singular points on the cylinder are shown in fig. 16.

\begin{figure}[ht!]
\center{\includegraphics[width=0.7\linewidth]{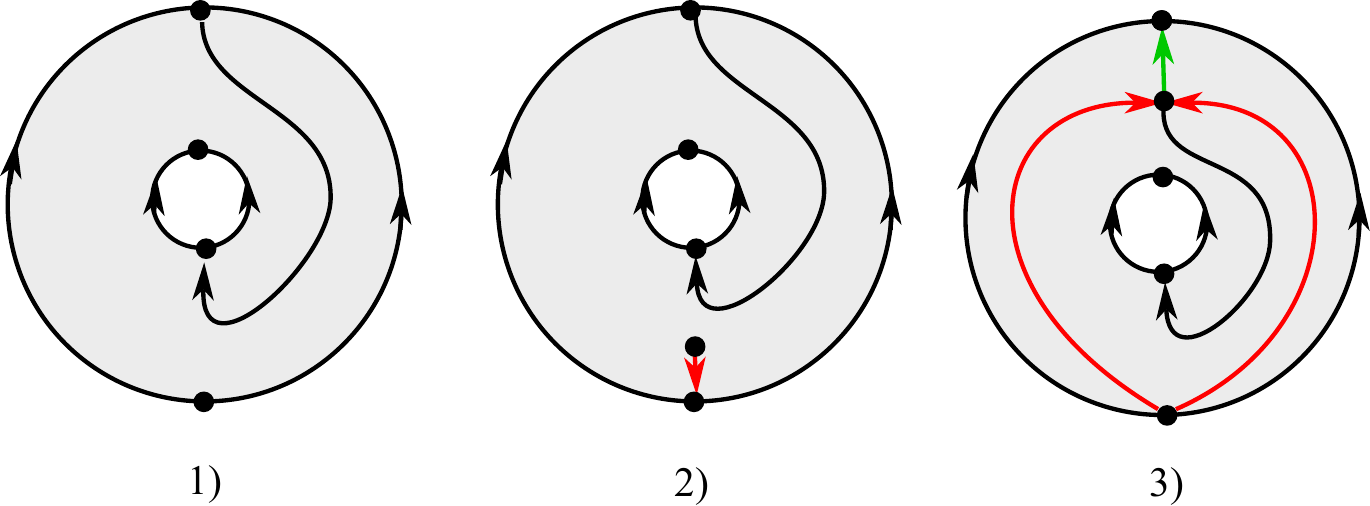}
}
{
\caption{saddle connections on a cylinder with 4 and 5 singular points}
}
\label{SI-sc45}
\end{figure}

Flows with 6 singular points and one saddle connection are shown in fig. 17.

\begin{figure}[ht!]
\center{\includegraphics[width=0.95\linewidth]{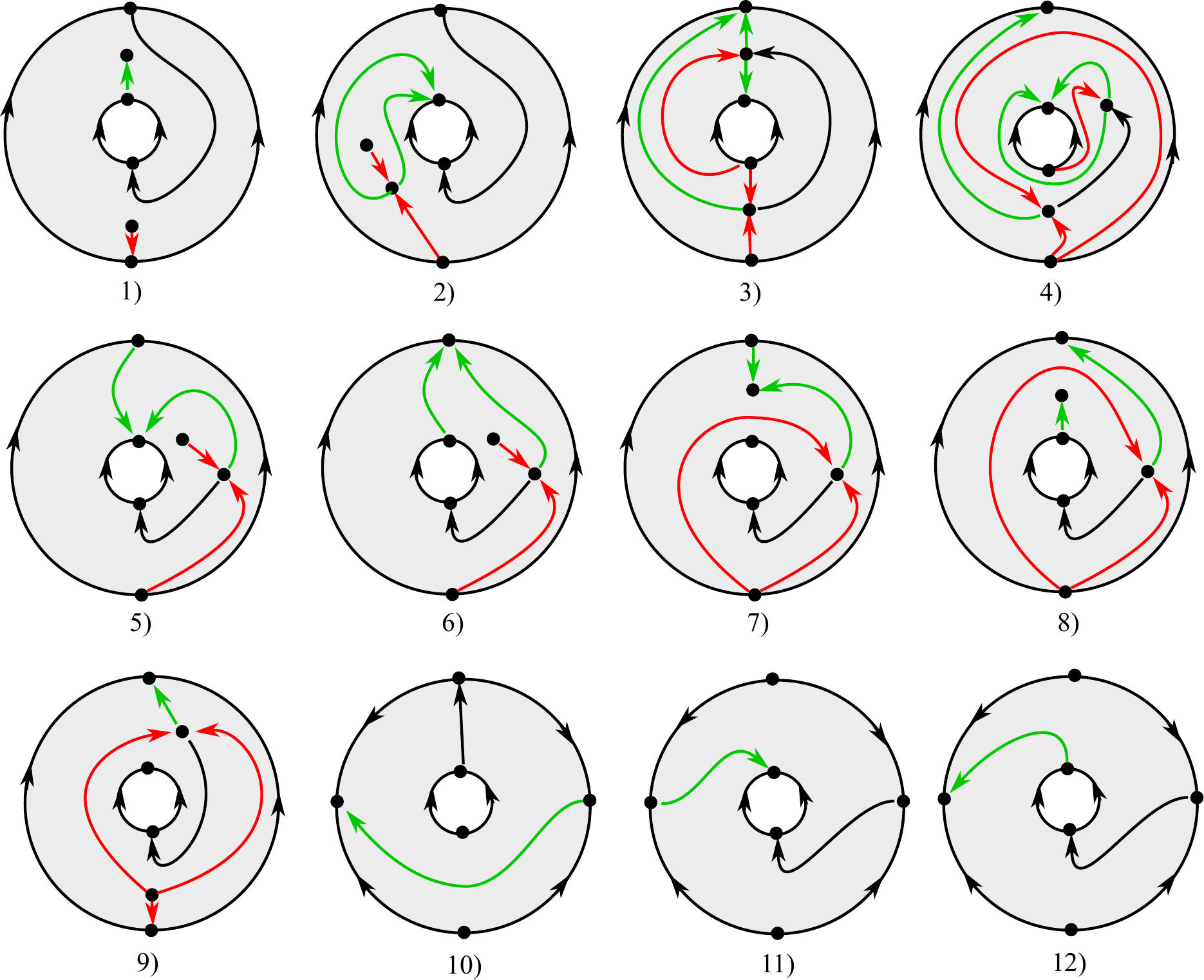}
}
{
\caption{saddle connections on a cylinder with 6 singular points}
}
\label{SI-sc6}
\end{figure}

We went through all the possible options, and therefore it is fair

\begin{theorem}
The following possible structures of typical one-parameter gradient saddle-node bifurcations exist on the cylinder:
\begin{itemize}
\item
with four singular points at the moment of bifurcation - ten bifurcations of the HN type;
\item
with five singular points at the moment of bifurcation - 4 SN, 4 HN, 18 HS, 14 BSN, 6 BDS.
\end{itemize}

There are the following bifurcations of the saddle connection:
\begin{itemize}
\item
with four singular points one BSC ( 1) in fig. 16);
\item
with five singular points -- two BSC ( 2) and reversed in fig. 16) and two HSCs ( 3) and inverted in fig. 16);
\item
with six singular points (Fig. 17) --
two SC (3),4)), ten HSC (5)$\times$2, 6)$\times$2, 7)$\times$2, 8)$\times$2, 9)$\times$2) and nine BSC (1,2) $\times$2,10) $\times$2, 11)$\times$2, 12) $\times$2).
\end{itemize}

\end{theorem}

\section{Structure of typical flows and bifurcations on a sphere with three holes}

A sphere with three holes will be considered as a three-connected region on the plane (a 2-disk with two holes).

Since there must be two singular points on each of the holes, streams with 6 singular points have no other singular points. Therefore, there will be no internal or semi-boundary bifurcations in this case.

The doubling of a sphere with three holes is a surface of genus 2. Except for one source and one sink, all other points are saddle points (otherwise the Poincaré-Hopf theorem on the sum of indices is violated).
There are only two possible structures of Morse streams in this case. In the first case, the source and the sink lie on one component of the boundary, and in the second - on different components (see Fig. 18).

\begin{figure}[ht!]
\center{\includegraphics[width=0.6\linewidth]{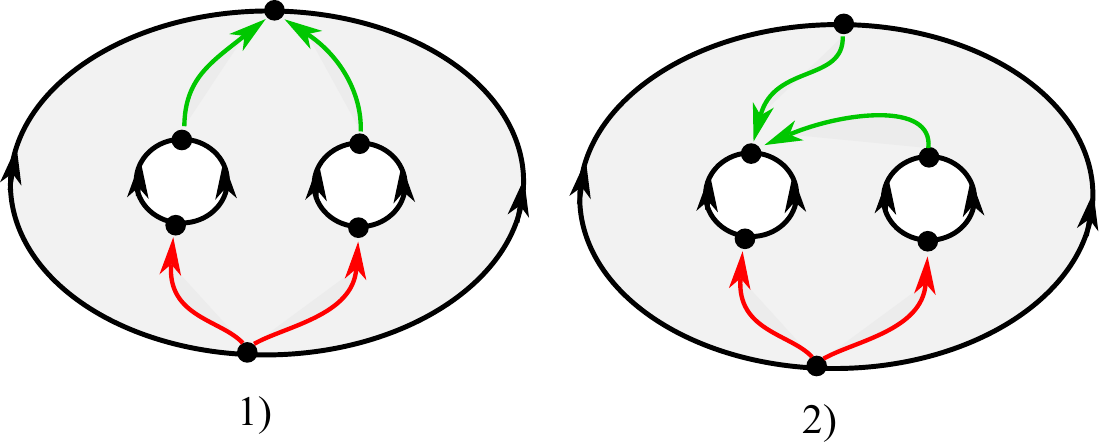}
}
{
\caption{Morse flows with 6 singular points on a sphere with three holes }
}
\label{F03}
\end{figure}

Since only two singular points lie on each boundary component and none of the separatrixes can be compressed without changing the topological type of the surface, saddle-node bifurcations do not exist in this case.

Bifurcations of the saddle connection are obtained when one separatrix slides over another.
All possible 4 cases are shown in fig. 19.

\begin{figure}[ht!]
\center{\includegraphics[width=0.6\linewidth]{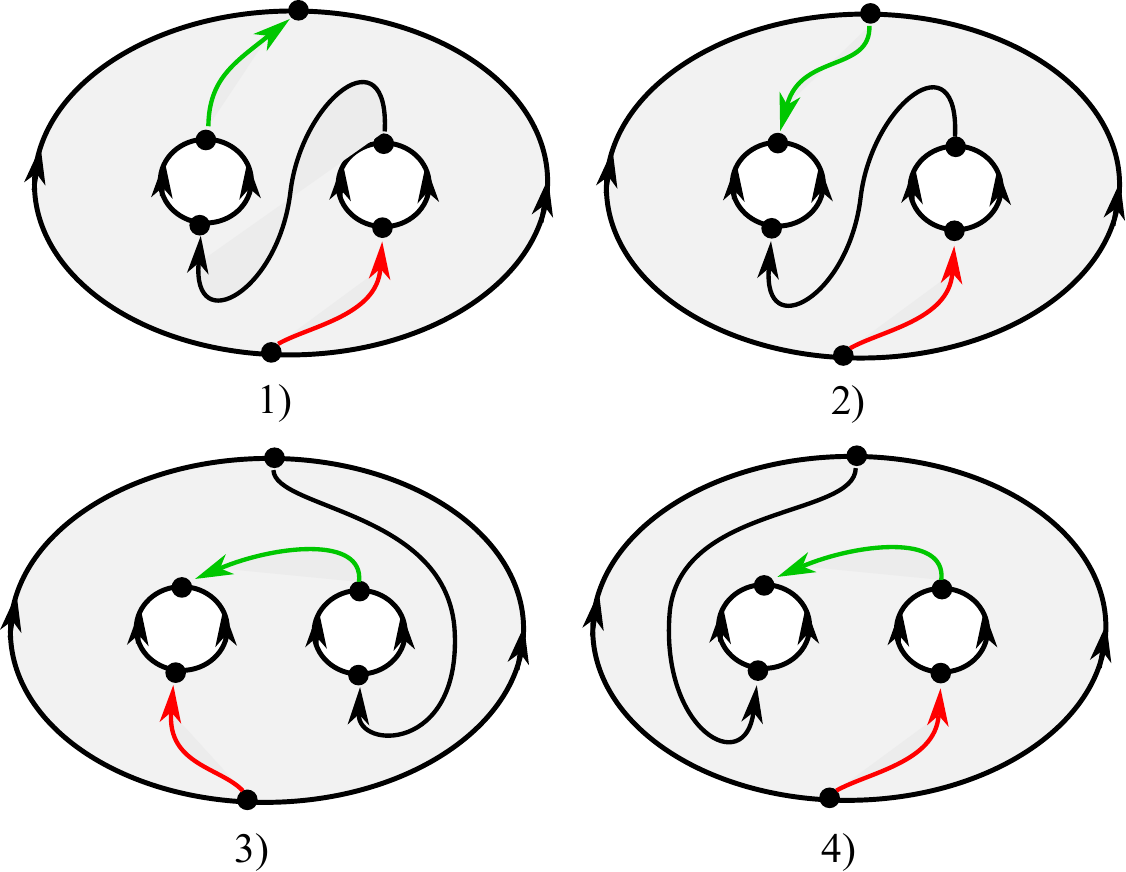}
}
{
\caption{saddle connections on a sphere with three holes }
}
\label{F03sc}
\end{figure}

Therefore, the following statement is true.

\begin{theorem}
On a sphere with three holes, there are no gradient saddle-node bifurcations with 5 singular points at the moment of bifurcation, but there are four boundary saddle bifurcations with 6 singular points.
\end{theorem}

\section*{Conclusion}

All possible structures of Morse flows and typical one-parameter bifurcations on spheres with holes in which no more than six singular points are found (see Table 1). We hope that the research carried out in this work can be extended to other surfaces and with a larger number of singular points.

\begin{table}[ht]
	\centering
		\begin{tabular} {|c|c|c|c|c|c|c|c|c|c|}
		\hline
Number of points
& 
Morse & SN & SC & BSN & BDS & HN & HS & HSC & BSC
\\ 
\hline
3 on $D^2$ & 2 & 0 & 0 & 0 & 0 & 2 & 0 & 0 & 0 
 \\
\hline
4 on $D^2$ & 5 & 2 & 0 & 2 & 0 & 0 & 2 & 4 & 0 
 \\
\hline
5 on $D^2$ & 7 & 8 & 0 & 2 & 0 & 6 & 8 & 4 & 0 
 \\
\hline
6 on $D^2$ & 22 & 30 & 7 & 22 & 5 & 12 & 38 & 6 & 2
 \\
\hline
4 on $S^1 \times I$ & 2 & 0 & 0 & 0 & 0 & 0 & 0 & 0 & 1
 \\
\hline
5 on $S^1 \times I$ & 4 & 0 & 0 & 0 & 10 & 0 & 0 & 2 & 2
 \\
\hline
6 on $S^1 \times I$& 14 & 4 & 2 & 14 & 6 & 4 & 18 & 10 & 9
 \\
\hline
6 on $F_{0,3}$ & 2 & 0 & 0 & 0 & 0 & 0 & 0 & 0 & 4

\\
 \hline		

		\end{tabular}
	\caption{Number of Morse flows and bifurcations on $S^2$ with holes (number of points befor bifurcation)}
	\label{tab:NF}
\end{table}


\begin{thebibliography}{10}

\bibitem{akchurin2022three}
O.~Akchurin, S.~Bilun, and A.~Prishlyak.
\newblock Three-color graph as the 1-skeleton of the 2-sphere triangulation.
\newblock {\em arXiv preprint arXiv:2209.05737}, 2022.
\href{http://dx.doi.org/10.48550/ARXIV.2209.05737 }{\path{doi: 10.48550/ARXIV.2209.05737}}.

\bibitem{bilun2002closed}
S.~Bilun and A.~Prishlyak.
\newblock The closed morse 1-forms on closed surfaces.
\newblock {\em Visn., Mat. Mekh., Kyv. Univ. Im. Tarasa Shevchenka},
  2002(8):77--81, 2002.

\bibitem{bilun2022visualization}
S.~Bilun and A.~Prishlyak.
\newblock Visualization of morse flow with two saddles on 3-sphere diagrams.
\newblock {\em arXiv preprint arXiv:2209.12174}, 2022.
\href{http://dx.doi.org/10.48550/ARXIV.2209.12174}{\path{doi: 10.48550/ARXIV.2209.12174}}.


\bibitem{bilun2022morse}
S.~Bilun, A.~Prishlyak, and A.~Prus.
\newblock Morse flows with fixed points on the boundary of 3-manifold.
\newblock {\em arXiv preprint arXiv:2209.04019}, 2022.
\href{http://dx.doi.org/10.48550/arXiv.2209.04019 }{\path{doi: 10.48550/arXiv.2209.04019}}.


\bibitem{bilun2023discrete}
S. Bilun, M. Hrechko, O. Myshnova, A. Prishlyak.
\newblock Structures of optimal discrete gradient vector fields on surface with one or two critical cells
\newblock {\em arXiv preprint arXiv:2303.07258}, 2023.
\href{http://dx.doi.org/10.48550/arXiv.2303.07258}{\path{doi: 10.48550/arXiv.2303.07258}}.


\bibitem{bilun2023morseRP2}
S. Bilun, A. Prishlyak, S. Stas, A. Vlasenko.
\newblock Topological structure of Morse functions on the projective plane
\newblock {\em arXiv preprint arXiv:2303.03850}, 2023.
\href{http://dx.doi.org/10.48550/arXiv.2303.03850}{\path{doi: 10.48550/arXiv.2303.03850}}.


\bibitem{bilun2023gradient}
S. Bilun, B. Hladysh, A. Prishlyak, V Sinitsyn.
\newblock Gradient vector fields of codimension one on the 2-sphere with at most ten singular points
\newblock {\em arXiv preprint arXiv:2303.10929}, 2023.
\href{http://dx.doi.org/10.48550/arXiv.2303.10929}{\path{doi: 10.48550/arXiv.2303.10929}}.



\bibitem{Bolsinov2004}
A.V. Bolsinov and A.T. Fomenko.
\newblock {\em Integrable Hamiltonian systems. Geometry, Topology,
  Classification}.
\newblock A CRC Press Company, Boca Raton London New York Washington, D.C.,
  2004.
\newblock 724 p.




\bibitem{hatamian2020heegaard}
C.~Hatamian and A.~Prishlyak.
\newblock Heegaard diagrams and optimal morse flows on non-orientable
  3-manifolds of genus 1 and genus 2.
\newblock {\em Proceedings of the International Geometry Center}, 13(3):33--48,
  2020.
\href{http://dx.doi.org/10.15673/tmgc.v13i3.1779}{\path{doi: 10.15673/tmgc.v13i3.1779}}.


\bibitem{hladysh2016functions}
B.~I. Hladysh and A.~O. Pryshlyak.
\newblock Functions with nondegenerate critical points on the boundary of the
  surface.
\newblock {\em Ukrainian Mathematical Journal}, 68(1):29--41, 2016.
\href{http://dx.doi.org/10.1007/s11253-016-1206-5}{\path{doi: 10.1007/s11253-016-1206-5}}.

\bibitem{hladysh2020deformations}
B.~I. Hladysh and A.~O. Pryshlyak.
\newblock Deformations in the general position of the optimal functions on oriented surfaces with boundary.
\newblock {\em Ukrainian Mathematical Journal}, 71(8):1173--1185, 2020.
\href{http://dx.doi.org/10.1007/s11253-019-01706-8}{\path{doi: 10.1007/s11253-019-01706-8}}.


\bibitem{hladysh2017topology}
B.I. Hladysh and A.O. Prishlyak.
\newblock Topology of functions with isolated critical points on the boundary
  of a 2-dimensional manifold.
\newblock {\em SIGMA. Symmetry, Integrability and Geometry: Methods and
  Applications}, 13:050, 2017.
\href{http://dx.doi.org/0.3842/SIGMA.2017.050}{\path{doi: 0.3842/SIGMA.2017.050}}.

\bibitem{hladysh2019simple}
B.I. Hladysh and A.O. Prishlyak.
\newblock Simple morse functions on an oriented surface with boundary.
\newblock {\em Журнал математической физики,
  анализа, геометрии}, 15(3):354--368, 2019.
\href{http://dx.doi.org/10.15407/mag15.03.354}{\path{doi: 10.15407/mag15.03.354}}.

\bibitem{Kronrod1950}
A.S. Kronrod.
\newblock Functions of two variables.
\newblock {\em Russian Mathematical Surveys}, 5:24--134, 1950.


\bibitem{Kybalko2018}
Z.~Kybalko, A.~Prishlyak, and R.~Shchurko.
\newblock {Trajectory equivalence of optimal Morse flows on closed surfaces}.
\newblock {\em {Proc. Int. Geom. Cent.}}, 11(1):12--26, 2018.
\href{http://dx.doi.org/10.15673/tmgc.v11i1.916 }{\path{doi: 10.15673/tmgc.v11i1.916}}.


\bibitem{loseva2016topology}
M.~Losieva and A.~Prishlyak.
\newblock Topology of morse--smale flows with singularities on the boundary of
  a two-dimensional disk.
\newblock {\em Pr. Mizhnar. Heometr. Tsentr}, 9(2):32--41, 2016.
\href{http://dx.doi.org/10.15673/tmgc.v9i2.279}{\path{doi: 10.15673/tmgc.v9i2.279}}.

\bibitem{lychak2009morse}
D.P. Lychak and A.O. Prishlyak.
\newblock Morse functions and flows on nonorientable surfaces.
\newblock {\em Methods of Functional Analysis and Topology}, 15(03):251--258,
  2009.



\bibitem{Oshemkov1998}
A.A. Oshemkov and V.V. Sharko.
\newblock Classication of morse-smale flows on two-dimensional manifolds.
\newblock {\em Matem. Sbornik}, 189(8):93--140, 1998.

\bibitem{Peixoto1973}
M.M. Peixoto.
\newblock On the classication of flows of 2-manifolds.
\newblock {\em Dynamical Systems (Proc. Symp. Univ. of Bahia, Salvador, Brasil,
  1971)}, pages 389--419, 1973.


\bibitem{prishlyak1997graphs}
A.O. Prishlyak.
\newblock On graphs embedded in a surface.
\newblock {\em Russian Mathematical Surveys}, 52(4):844, 1997.
\href{http://dx.doi.org/10.1070/RM1997v052n04ABEH002074}{\path{doi: 10.1070/RM1997v052n04ABEH002074}}.

\bibitem{prishlyak2001conjugacy}
A.O. Prishlyak.
\newblock Conjugacy of Morse functions on 4-manifolds.
\newblock {\em Russian Mathematical Surveys}, 56(1):170, 2001.




\bibitem{prishlyak2002ms}
A.O. Prishlyak.
\newblock Morse--smale vector fields without closed trajectories on-manifolds.
\newblock {\em Mathematical Notes}, 71(1-2):230--235, 2002.
\href{http://dx.doi.org/10.1023/A:1013963315626}{\path{doi: 10.1023/A:1013963315626}}.

\bibitem{prishlyak2003sum}
A.O. Prishlyak.
\newblock On sum of indices of flow with isolated fixed points on a stratified
  sets.
\newblock {\em Zhurnal Matematicheskoi Fiziki, Analiza, Geometrii [Journal of
  Mathematical Physics, Analysis, Geometry]}, 10(1):106--115, 2003.

\bibitem{prishlyak2005complete}
A.O. Prishlyak.
\newblock Complete topological invariants of morse--smale flows and handle
  decompositions of 3-manifolds.
\newblock {\em Fundamentalnaya i Prikladnaya Matematika}, 11(4):185--196, 2005.

\bibitem{prishlyak2007complete}
A.O. Prishlyak.
\newblock Complete topological invariants of morse-smale flows and handle
  decompositions of 3-manifolds.
\newblock {\em Journal of Mathematical Sciences}, 144:4492--4499, 2007.

\bibitem{prishlyak2022Boy}
A.~Prishlyak and L.~Di Beo.
\newblock Flows with minimal number of singularities on the Boy’s surface
\newblock {\em Proceedings of the International Geometry Center}, 15(1):32--49,
  2020.

\bibitem{prishlyak2022topological}
A.~Prishlyak and M.~Loseva.
\newblock Topological structure of optimal flows on the girl's surface.
\newblock {\em Proceedings of the International Geometry Center},
  15(3-4):184--202, 2022.

\bibitem{prishlyak2021flows}
A.~Prishlyak, A.~Prus, and S.~Huraka.
\newblock Flows with collective dynamics on a sphere.
\newblock {\em Proc. Int. Geom. Cent}, 14(1):61--80, 2021.
\href{http://dx.doi.org/10.15673/tmgc.v14i1.1902}{\path{doi: 10.15673/tmgc.v14i1.1902}}.




\bibitem{prishlyak2020topology}
A.~Prishlyak and M.~Loseva.
\newblock Topology of optimal flows with collective dynamics on closed
  orientable surfaces.
\newblock {\em Proceedings of the International Geometry Center}, 13(2):50--67,
  2020.
\href{http://dx.doi.org/10.15673/tmgc.v13i2.1731}{\path{doi: 10.15673/tmgc.v13i2.1731}}.

\bibitem{prishlyak2017morse}
A.~Prishlyak and A.~Prus.
\newblock Morse-smale flows on torus with hole.
\newblock {\em Proc. Int. Geom. Cent.}, 10(1):47--58, 2017.
\href{http://dx.doi.org/10.15673/tmgc.v1i10.549}{\path{doi: 10.15673/tmgc.v1i10.549}}.

\bibitem{prishlyak1999equivalence}
A.O.~Prishlyak.
\newblock Equivalence of morse function on 3-manifolds.
\newblock {\em Methods of Func. Ann. and Topology}, 5(3):49--53, 1999.

\bibitem{prishlyak2000conjugacy}
A.O.~Prishlyak.
\newblock Conjugacy of morse functions on surfaces with values on a straight
  line and circle.
\newblock {\em Ukrainian Mathematical Journal}, 52(10):1623--1627, 2000.
\href{http://dx.doi.org/10.1023/A:1010461319703}{\path{doi: 10.1023/A:1010461319703}}.

\bibitem{prishlyak2002morse}
A.O.~Prishlyak.
\newblock Morse functions with finite number of singularities on a plane.
\newblock {\em Meth. Funct. Anal. Topol}, 8:75--78, 2002.

\bibitem{Prishlyak2002beh2}
A.O.~Prishlyak.
\newblock Topological equivalence of morse--smale vector fields with beh2 on
  three-dimensional manifolds.
\newblock {\em Ukrainian Mathematical Journal}, 54(4):603--612, 2002.

\bibitem{prishlyak2003topological}
A.O.~Prishlyak.
\newblock Topological classification of m-fields on two-and three-dimensional
  manifolds with boundary.
\newblock {\em Ukrainian Mathematical Journal}, 55(6):966--973, 2003.

\bibitem{prishlyak2019optimal}
A.O.~Prishlyak and M.V.~Loseva.
\newblock Optimal morse--smale flows with singularities on the boundary of a
  surface.
\newblock {\em Journal of Mathematical Sciences}, 243:279--286, 2019.

\bibitem{prishlyak2007classification}
A.O.~Prishlyak and K.I.~Mischenko.
\newblock Classification of noncompact surfaces with boundary.
\newblock {\em Methods of Functional Analysis and Topology}, 13(01):62--66,
  2007.

\bibitem{prishlyak2020three}
A.O.~Prishlyak and A.A.~Prus.
\newblock Three-color graph of the morse flow on a compact surface with
  boundary.
\newblock {\em Journal of Mathematical Sciences}, 249(4):661--672, 2020.
\href{http://dx.doi.org/10.1007/s10958-020-04964-1}{\path{doi: 10.1007/s10958-020-04964-1}}.

\bibitem{Reeb1946}
G.~Reeb.
\newblock Sur les points singuliers d’une forme de pfaff complétement
  intégrable ou d’une fonction numérique.
\newblock {\em C.R.A.S. Paris}, 222:847—849, 1946.

\bibitem{Sharko1993}
V.V. Sharko.
\newblock {\em Functions on manifolds. Algebraic and topological aspects.},
  volume 131 of {\em Translations of Mathematical Monographs}.
\newblock American Mathematical Society, Providence, RI, 1993.

\bibitem{Smale1961}
S.~Smale.
\newblock On gradient dynamical systems.
\newblock {\em Ann. of Math.}, 74:199--206, 1961.

\bibitem{kkp2013}
В.М. Кузаконь, В.Ф. Кириченко, and О.О. Пришляк.
\newblock Гладкi многовиди. Геометричнi та
  топологiчнi аспекти.
\newblock {\em Працi Iнcтитуту математики НАН
  України.—2013.—97.—500 с}, 2013.

\bibitem{prish1998vek}
А.О. Пришляк.
\newblock Векторные поля Морса--Смейла с
  конечным числом особых траекторий на
  трехмерных многообразиях.
\newblock {\em Доповiдi НАН України}, (6):43--47, 1998.

\bibitem{prish1998sopr}
А.О. Пришляк.
\newblock Сопряженность функций Морса.
\newblock {\em Некоторые вопросы совр.
  математики. Институт математики АН
  Украины, Киев}, 1998.

\bibitem{prish2001top}
А.О. Пришляк.
\newblock Топологическая эквивалентность
  функций и векторных полей Морса—Смейла
  на трёхмерных многообразиях.
\newblock {\em Топология и геометрия. Труды
  Украинского мат. конгресса}, pages 29--38, 2001.

\bibitem{prish2002vek}
А.О. Пришляк.
\newblock Векторные поля Морса--Смейла без
  замкнутых траекторий на трехмерных
  многообразиях.
\newblock {\em Математические заметки}, 71(2):254--260,
  2002.

\bibitem{prish2015top}
О.О. Пришляк.
\newblock Топологія многовидів.
\newblock {\em Київський університет}, 2015.

\end{thebibliography}

\textsc{Taras Shevchenko National University of Kyiv}

Svitlana Bilun  \ \ \ \ \ \ \ \ \ \ \ \
\textit{Email address:} \text{ bilun@knu.ua}   \ \ \ \ \ \ \ \ \
\textit{ Orcid ID:}  \text{0000-0003-2925-5392}

Maria Loseva \ \ \ \ \ \ \ 
\textit{Email:} \text{ mv.loseva@gmail.com} \   
\textit{ Orcid ID:} \text{0000-0002-2282-206X}

Olena Myshnova \ \ \ \ \ \ \ \ \
\textit{Email:} \text{ myshnova.olena@gmail.com}  \ 
\textit{ Orcid ID:} \text{0009-0001-8808-0038}

Alexandr Prishlyak \ \ \ \ \ 
\textit{Email address:} \text{ prishlyak@knu.ua} \ \ \ \ 
\textit{ Orcid ID:} \text{0000-0002-7164-807X}

\end{document}